\documentstyle[11pt,twoside]{article}

\newtheorem{theorem}{\noindent\bf Theorem}[section]
\newtheorem{proposition}[theorem]{\noindent\bf Proposition}
\newtheorem{lemma}[theorem]{\noindent\bf Lemma}
\newtheorem{corollary}[theorem]{\noindent\bf Corollary}

\newtheorem{example}[theorem]{\noindent\bf Example}

\newcommand{\supp}{\hbox{supp}}

\newcommand{\e}{\hfill\blacksquare}
\input{amssym}
\begin{document}
\date{}
\title{{\Large\bf On the weak$^*$ continuity of $LUC({\cal G})^*$-module
action on $LUC({\cal X},{\cal G})^*$ related to ${\cal G}$-space $\cal
X$} }
\author{{\normalsize\sc H. Javanshiri and N. Tavallaei}}

\maketitle
\normalsize

\begin{abstract}
Associated with a locally compact group $\cal G$ and a $\cal G$-space
$\cal X$ there is a Banach subspace $LUC({\cal X},{\cal G})$ of $C_b({\cal X})$,
which has been introduced and studied by Lau and Chu in \cite{chulau}.
In this paper,
we study some properties of the first dual space of $LUC({\cal
X},{\cal G})$. In particular, we introduce a left action of
$LUC({\cal G})^*$ on $LUC({\cal X},{\cal G})^*$ to make it a
Banach left module and then we investigate the Banach subalgebra
${{\frak{Z}}({\cal X},{\cal G})}$ of $LUC({\cal G})^*$, as the
topological centre related to this module action, which contains
$M({\cal G})$ as a closed
subalgebra. Also, we show that the faithfulness of this module action
is related to the properties of the action of $\cal G$ on $\cal X$ and we
extend the main results of Lau~\cite{lau} from locally compact
groups to ${\cal G}$-spaces. Sufficient and/or necessary conditions for the equality
${{\frak{Z}}({\cal X},{\cal G})}=M({\cal G})$ or $LUC({\cal G})^*$ are given.
Finally, we apply our
results to some special cases of $\cal G$ and $\cal X$ for obtaining
various examples whose topological centres ${{\frak{Z}}({\cal X},{\cal G})}$
are $M({\cal G})$, $LUC({\cal G})^*$ or neither of them.\\


\end{abstract}

%%%%%%%%%%%%%%%%%%%%%%%%%%%%%%%%%%%%%%%%%%%%%%%%%%%%%%%%%%%%%%%%%%%%%%%%%%%%%%%%%%%%%%%%%%%%%%%%%%
%%%%%%%%%%%%%%%%%%%%%%%%%%%%%%%%%%%%%%%%%%%%%%%%%%%%%%%%%%%%%%%%%%%%%%%%%%%%%%%%%%%%%%%%%%%%%%%%%%

\section{Introduction}

Let $\cal G$ be a locally compact group. Then the Banach space
$LUC({\cal G})^*$, the topological dual of the space of all bounded
left uniformly continuous functions on $\cal G$, is a Banach
algebra equipped with the Arens-type product. In fact, the
product of $m$ and $n$ in $LUC({\cal G})^*$ is defined by
\begin{eqnarray}\label{arensluc}
\langle m\odot n,f\rangle =\langle m,n\diamond f\rangle ,
\quad\quad (n\diamond f)(s)=\langle n,{_sf}\rangle ,\quad\quad {_sf}(t)=f(st),
\end{eqnarray}
where $f\in LUC({\cal G})$ and $s,t\in {\cal G}$. In general, this
product is not separately weak$^*$ to weak$^*$ continuous on
$LUC({\cal G})^*$, and in recent years there has been shown
considerable interest by harmonic analysts in the characterization
of the following space
$${\frak{Z}}({\cal G})=\Big{\{}m\in LUC({\cal G})^*:~n\mapsto m\odot n~
{\hbox{is}}~{\hbox{weak}}^*~{\hbox{to}}~{\hbox{weak}}^*
~{\hbox{continuous}}\Big{\}}.$$ As far as we know the subject, the
starting point of the study of the space ${\frak{Z}}({\cal G})$ is
the paper by Zappa \cite{zappa}. In details, Zappa proved that
${\frak{Z}}({\Bbb R})$ is precisely $M({\Bbb R})$, where $\Bbb R$
is the additive group of real number and $M({\Bbb R})$ is the
Banach space of all complex Radon measures on $\Bbb R$. This
result was extended to all abelian locally campact groups by
Grosser and Losert in \cite{grossert}, and to all locally compact
groups by Lau in \cite{lau}.

%%%%%%%%%%%%%%%%%%%%%%%%%%%%%%%%%%%%%%%%%%%%%%%%%%%%%%%%%%%%%%%%%%%%%%%%%%%%%%%%%%%%%%%%%%%%%%%

In this paper, considering $\cal X$ as a locally compact Hausdorff
space on which $\cal G$ acts continuously from the left, following Lau and Chu \cite{chulau} we
introduce the Banach space $LUC({\cal X},{\cal G})$.
Then, we present a left action of $LUC({\cal
G})^*$ on $LUC({\cal X},{\cal G})^*$, as an extension of the
natural action of $M({\cal G})$ on $M({\cal X})$, to make a Banach
left $LUC({\cal G})^*$-module. In
particular, we investigate the faithfulness of the action of $M({\cal
G})$ (resp. $LUC({\cal G})^*$) on $M({\cal X})$ (resp. $LUC({\cal X},{\cal G})^*$)
in relation to the action of $\cal G$ on $\cal X$.
Also, we prove that if $\cal X$ is a transitive $\cal
G$-space, then $M({\cal X})$ is an $LUC({\cal G})^*$-submodule of
$LUC({\cal X},{\cal G})^*$ just when $\cal X$ is compact.
Moreover, the main purpose of this work is to study the
topological centre problem related to this module action. To this end,
we introduce the Banach subspace ${{\frak{Z}}({\cal X},{\cal G})}$
of $LUC({\cal G})^*$, as the topological centre related to this
module action, which contains $M({\cal G})$ as a closed subalgebra.
Furthermore, apart from some characterization of the space
${{\frak{Z}}({\cal X},{\cal G})}$, we extend the main results of
Lau~\cite{lau} from locally compact groups to $\cal G$-spaces.  Finally,
we apply our results  to some
special cases of $\cal G$ and $\cal X$ to give examples with

$\centerdot$ ${\frak Z}({\cal X},{\cal G})=M({\cal G})$;

$\centerdot$ $M({\cal G})\subsetneqq{\frak Z}({\cal X},{\cal G})=LUC({\cal G})^*$;

$\centerdot$ $M({\cal G})\subsetneqq{\frak Z}({\cal X},{\cal G})\subsetneqq LUC({\cal G})^*$.

%%%%%%%%%%%%%%%%%%%%%%%%%%%%%%%%%%%%%%%%%%%%%%%%%%%%%%%%%%%%%%%%%%%%%%%%%%%%%%%%%%%%%%%%%%%%%%%

The motivation of the study presented in this paper comes from
some works of Lau and his
coauthors \cite{chulau,lauucb1,lauucb2,lau,laupym1} and some recent developments of the concept of
regularity of bilinear maps were introduced and first studied by
Arens \cite{arens}. It is known that the Arens regularity of module
actions has been a major tool in the study of Banach algebras. For
example, Arens regularity of module actions as a generalization of
Arens regularity of Banach algebras were introduced and first
studied by Filali and Eshaghi Gordji \cite{filalgordji} and they
used this notion to answer some questions regarding Arens
regularity of Banach algebras raised by Lau and \"{U}lger
\cite{laulger}. Also in \cite{mvishki}, Arens regularity of module
actions were considered by Mohammadzadeh and Vishki to investigate
the conditions under which the second adjoint of a derivation into
a dual Banach module is again a derivation, that extends the
results of Dales, Rodriguez-Palacios and Velasco in
\cite{dalesrpv} for a general derivation.

%%%%%%%%%%%%%%%%%%%%%%%%%%%%%%%%%%%%%%%%%%%%%%%%%%%%%%%%%%%%%%%%%%%%%%%%%%%%%%%%%%
%%%%%%%%%%%%%%%%%%%%%%%%%%%%%%%%%%%%%%%%%%%%%%%%%%%%%%%%%%%%%%%%%%%%%%%%%%%%%%%%%%

\section{\normalsize\bf Some prerequisites}

Throughout this paper, $\cal G$ is a locally compact group with
left Haar measure $\lambda$ and identity element $e$,
$\Delta$ refers to the modular function on $\cal G$ and
the notation $C_b({\cal G})$ is used to denote the space of all
bounded complex-valued continuous functions on $\cal G$ with the
supremum norm. Also, the subspaces $C_c({\cal G})$ and $C_0({\cal
G})$ of $C_b({\cal G})$ refer to the space of all functions with
compact support and the space of all functions vanishing at infinity,
respectively. Also, for all $s\in {\cal G}$ and $f\in C_b({\cal
G})$ we define ${_sf}, $ the left translation of $f$ by $s$, as
${_sf}(t)=f(st)$, $t\in {\cal G}$. A function $f$ in $C_b({\cal
G})$ is called left uniformly continuous if the mapping
$s\mapsto{_sf}$ from $\cal G$ into $(C_b({\cal
G}),\|\cdot\|_\infty)$ is continuous. As usual, we mean by
$LUC({\cal G})$ the Banach space of all left uniformly continuous
functions on $\cal G$. Moreover, the notation $M({\cal G})$
is used to denote the measure algebra of ${\cal G}$ consisting of all complex regular
Borel measures on ${\cal G}$ with the total variation norm and the
convolution product ``$\ast$" defined by the formula
 $$
 \langle\mu\ast\nu ,f\rangle=\int_{\cal G}\int_{\cal G} f(st)\;
 d\mu(s)\; d\nu(t)
 $$
for all $\mu, \nu\in M({\cal G})$ and $f\in C_0({\cal G})$.
It is well-known that $M({\cal G})$
is the topological dual of $C_0({\cal G})$ with the pairing
 $$
 \langle \mu,f \rangle:=\int_{\cal G}  f(s)\;d\mu(s)
 $$
for all $\mu\in M({\cal G})$ and $f\in C_0({\cal G})$.
This allows us to view $M({\cal G})$ as a subalgebra of $LUC({\cal
G})^*$.

%%%%%%%%%%%%%%%%%%%%%%%%%%%%%%%%%%%%%%%%%%%%%%%%%%%%%%%%%%%%%%%%%%%%%%%%%%%%%%%%%%%%%%%%%%%%%%%

Throughout the paper the letter $M^+({\cal G})$ means the set of
all positive measures in $M({\cal G})$; $\delta_s$ denotes
the Dirac measure at $s\in{\cal G}$; and as usual, $M_a({\cal G})$ denotes the closed ideal of $M({\cal
G})$ consisting of all absolutely continuous measures with respect
to $\lambda$. Let also $L^1({\cal G})$ denote the group algebra of
${\cal G}$ as defined in \cite{farmonic}.
Then, the Radon-Nikodym theorem can be interpreted as an
identification  of $M_a({\cal G})$ with
${\Big\{}\lambda_\varphi:~\varphi\in L^1({\cal G}){\Big\}}$, where
$\lambda_\varphi$ is the measure in $M({\cal G})$ defined on each
Borel subset $B$ of ${\cal G}$ by
$$\lambda_\varphi(B)=\int_B\varphi\;d\lambda.$$
Recall that, the convolution of two functions $\varphi, \psi\in
L^1({\cal G})$ is the function defined by
\begin{eqnarray*}
\varphi\ast\psi(s)&=&\int_{\cal
G}\varphi(t)\psi(t^{-1}s)\;d\lambda(t)\\
&=&\int_{\cal G}\varphi(st^{-1})\psi(t)\Delta(t^{-1})\;d\lambda(t).
\end{eqnarray*}

%%%%%%%%%%%%%%%%%%%%%%%%%%%%%%%%%%%%%%%%%%%%%%%%%%%%%%%%%%%%%%%%%%%%%%%%%%%%%%%%%%
%%%%%%%%%%%%%%%%%%%%%%%%%%%%%%%%%%%%%%%%%%%%%%%%%%%%%%%%%%%%%%%%%%%%%%%%%%%%%%%%%%

\section{\normalsize\bf Definitions and some basic results}

Recall that a locally compact Hausdorff space $\cal X$ is said to be a ({\it
left}) $\cal G$-{\it space} if there is a continuous action map
$${\cal G}\times {\cal X}\rightarrow {\cal X};~(s,x)\mapsto s\cdot x,$$
satisfying $(st)\cdot x=s\cdot(t\cdot x)$
and $e\cdot x=x$, for all $s,t\in {\cal G}$ and $x\in {\cal X}$.
In some published articles the terms ({\it left}) {\it
transformation group}, a {\it dynamical system} or a {\it flow}
are also used for the pair $({\cal X},{\cal G})$. Similarly one
can define a right action and a right $\cal G$-space. Through this
paper, when we refer to $\cal X$ as a $\cal G$-space without any explicit
reference to the left or right action, we mean $\cal X$ is a left
$\cal G$-space.

A $\cal G$-space $\cal X$ is called transitive if for each $x,y\in
{\cal X}$ there exists an element $s\in {\cal G}$ such that
$y=s\cdot x$. Also, in the case where $\cal X$ is a $\cal
G$-space which is topologically isomorphic to ${\cal G}/{\cal H}$,
for some closed subgroup $\cal H$ of $\cal G$, we say $\cal
X$ is a homogeneous $\cal G$-space.

%%%%%%%%%%%%%%%%%%%%%%%%%%%%%%%%%%%%%%%%%%%%%%%%%%%%%%%%%%%%%%%%%%%%%%%%%%%%%%%%%%%%%%%%%%%%%%%

The action of $\cal G$ on $\cal X$ is said to be {\it free}, where
$s\cdot x\neq x$ for every $x\in{\cal X}$ and $s\in{\cal G}$ with
$s\neq e$. It is called {\it effective}, when the following closed normal
subgroup of $\cal G$ is trivial
\[{\frak N}({\cal X},{\cal G})={\Big\{}s\in {\cal G}:\,
 s\cdot x=x{\hbox{~for~all}}\; x\in {\cal X}{\Big\}}.\]

%%%%%%%%%%%%%%%%%%%%%%%%%%%%%%%%%%%%%%%%%%%%%%%%%%%%%%%%%%%%%%%%%%%%%%%%%%%%%%%%%%%%%%%%%%%%%%%

Similar to the previous section, we consider the Banach spaces
$C_b({\cal X})$ and $C_0({\cal X})$. We also denote by $M({\cal X})$ the Banach space of all complex Radon measures on $\cal X$ with total variation norm.
For given measures
$\mu\in M({\cal G})$ and $\sigma\in M({\cal X})$ we may define
 their convolution $\mu\star\sigma$, as an element of $M({\cal X})$, by
\begin{eqnarray*}
\langle\mu\star\sigma,F\rangle=\int_{\cal G}\int_{\cal X} F(s\cdot
x)\;d\sigma(x)\;d\mu(s)\quad\quad\quad(F\in C_0({\cal X})).
\end{eqnarray*}
This allows us to consider $M({\cal X})$ as a Banach left-$M({\cal G})$-module. Now let $\cal A$ be a closed subalgebra of $M({\cal G})$, then as usual, we say that $M({\cal X})$ is a faithful Banach left
$\cal A$-module, if the action of each $0\neq\mu\in{\cal A}$ is
non-trivial; that is, if $\mu\in {\cal A}$ is so that $\mu\star
\sigma=0$ for all $\sigma\in M({\cal X})$, then $\mu=0$.
In the following result, we are concerned with relations between
faithfulness of the action of $M({\cal G})$ on $M({\cal X})$ and
the action of $\cal G$ on $\cal X$. In a special case, it gives some
conditions in relation to Veech's Theorem \cite{veech} in terms of the
faithfulness of the action of $M({\cal G})$ on $M({\cal U}{\cal
G})$ which is a functional property, where ${\cal U}{\cal G}$ is
the largest semigroup compactification of $\cal G$.

Let $(f_\alpha)_\alpha$ be a net in $C_b({\cal G})$. We say $(f_\alpha)_\alpha$ converges to some $f\in C_b({\cal G})$ strictly if $\|f_\alpha g-fg\|_\infty\longrightarrow 0$, for all $g$ in $S_0^+({\cal G})$ the set of all non-negative upper semicontinuous real-valued functions on $\cal G$ which vanish at infinity, see \cite{abel} for more details.

\begin{theorem}\label{210422}
Let $\cal G$ be a locally compact group and $\cal X$ be a $\cal
G$-space. Then the following statements hold.

\newcounter{j1212}
\begin{list}%
{\bf(\alph{j1212})}{\usecounter{j1212}}
\item The Banach space $M({\cal X})$ is a faithful Banach left $M^+({\cal
G})$-module if and only if $\cal G$ acts effectively on $\cal X$.
\item If for some $x_0\in {\cal G}$ the stabilizer subgroup
${\frak N}(x_0,{\cal G})=\{s\in {\cal G}:\,s\cdot x_0=x_0\}$ of ${\cal G}$
is trivial, then $M({\cal X})$ is a faithful Banach left $M({\cal
G})$-module. In particular, if  $\cal G$ acts freely on $\cal X$,
then $M({\cal G})$ acts faithfully on $M({\cal X})$.
\item If $M_a({\cal G})$ acts faithfully
on $M({\cal X})$, then $M({\cal X})$ is a faithful Banach left $M({\cal
G})$-module.
\end{list}
\end{theorem}

{\noindent Proof.} {\bf (a)} First, suppose that $M({\cal X})$ is a faithful
Banach left $M^+({\cal G})$-module and $s\in {\frak N}({\cal
X},{\cal G})$. Then for all $\sigma \in M({\cal X})$ we have
$$
\langle \delta_{s}\star \sigma , F\rangle=\int_{\cal X} F(s\cdot
x)\,d\sigma(x) =\int_{\cal X} F( x)\,d\sigma(x)= \langle \sigma ,
F\rangle,\quad\quad\quad {\Big(}F\in C_0({\cal X}){\Big)}
$$
Accordingly, $\delta_{s}\star \sigma=\sigma$ for all $\sigma \in
M({\cal X})$. This together with the fact that $M^+({\cal G})$ acts
faithfully on $M({\cal X})$ implies that $\delta_s=\delta_e$ and therefore $s=e$.

%%%%%%%%%%%%%%%%%%%%%%%%%%%%%%%%%%%%%%%%%%%%%%%%%%%%%%%%%%%%%%%%%%%%%%%%%%%%%%%%%%%

Conversely, assume that ${\frak N}({\cal X},{\cal G})=\{e\}$ and
$\mu$ is a positive measure in $M({\cal G})$ such that $\mu\star\sigma=0$ for all
$\sigma\in M({\cal X})$.
Observe that $\nu=\mu+\delta_e$ is a positive measure in $M({\cal G})$
for which $\nu\star \sigma=\sigma$ for all $\sigma\in
M({\cal X})$. Hence, the proof will be completed by showing that
$\supp(\nu)=\{e\}$ or equivalently
$\nu=\delta_e$. To this end, suppose on the contrary that $s_0$ is an element in $\supp(\nu)\setminus\{e\}$. This, together with the fact that
${\frak N}({\cal X},{\cal G})=\{e\}$ implies that
there exists some $x_0\in {\cal X}$ such that $s_0\cdot x_0\neq
x_0$. Taking a positive function $F$ in $C_c({\cal X})$ with
$F(x_0)=0$ and $F(s_0\cdot x_0)=1$, we have $$(\nu\star
\delta_{x_0})(F)=\delta_{x_0}(F)=F(x_0)=0.$$
On the
other hand, there is an open neighborhood $U$ of $s_0$ such that
$F(s\cdot x_0)>1/2$ for all $s\in U$. As $s_0\in \supp(\nu)$, we
have $\nu(U)>0$ and hence, we can write
\begin{eqnarray*}
0 = (\nu\star \delta_{x_0})(F)= \int_{\cal G} F(s\cdot
x_0)\,d\nu(s)\geq \nu(U)/2,
\end{eqnarray*}
which is a contradiction. Thus, we get $\supp(\nu)=\{e\}$. It follows that
$\nu=\delta_e$ and therefore $\mu=0$.

%%%%%%%%%%%%%%%%%%%%%%%%%%%%%%%%%%%%%%%%%%%%%%%%%%%%%%%%%%%%%%%%%%%%%%%%%%%%%%%%%%%%%

{\bf (b)} Assume that ${\frak N}(x_0,{\cal G})=\langle e\rangle$, for some $x_0\in {\cal X}$.
Obviously, for all $F\in C_b({\cal X})$, the function $r_{x_0}F:{\cal G}\rightarrow \Bbb{C}$ defined by $r_{x_0}F(s)=F(s\cdot x_0)$ belongs to $C_b({\cal G})$. Moreover,
$${\cal A}:={\Big\{}r_{x_0}F:\,F\in C_b({\cal X}){\Big\}},$$ is a self-adjoint subalgebra of $C_b({\cal G})$ such that for each $s\in {\cal G}$, there exists some $F\in C_b({\cal X})$ for which $r_{x_0}F(s)=F(s\cdot x_0)\neq 0$. Furthermore, this subalgebra separates the points of ${\cal G}$, since for all $s_1$ and $s_2$ in ${\cal G}$ with $s_1\neq s_2$ we have $s_1\cdot x_0\neq s_2\cdot x_0$ and therefore there exists $F\in C_b({\cal X})$ such that $$r_{x_0}F(s_1)=F(s_1\cdot x_0)\neq F(s_2\cdot x_0)=r_{x_0}F(s_2).$$  From this, by a generalized version of the Stone-Weierstrass theorem \cite[Theorem~3.2]{abel}, we can conclude that ${\cal A}$ is strictly dense in $C_b({\cal G})$.

Now, let $\mu$ be a measure in $M({\cal G})$ such that
$\mu\star \sigma=0$ for all $\sigma\in M({\cal X})$ and
$g\in C_c({\cal G})$ be given. Then, by taking $\sigma=\delta_{x_0}$, we have
\begin{eqnarray}\label{rm2}
\int_{\cal G} F(s\cdot x_0)\,d\mu(s) =  \langle\mu\star\delta_{x_0} , F\rangle=0,
\end{eqnarray}
for all $F\in C_b({\cal X})$.
Take a neighborhood $V$ of $\supp(g)$ with compact closure and some $F'\in C_c^+({\cal X})$ such that
$$F'{\Big|}_{\supp(g)\cdot x_0}=1\quad\quad{\hbox{and}}\quad\quad F'{\Big|}_{{\cal X}\setminus\overline{V}\cdot x_0}=0.$$ Therefore, $r_{x_0}F'$ is a non-negative continuous function for which
$$(r_{x_0}F'){\Big|}_{\supp(g)}=1\quad\quad{\hbox{and}}\quad\quad(r_{x_0}F'){\Big|}_{{\cal X}\setminus\overline{V}}=0;$$ This is because of, for all $s\in \cal G$, the condition $s\notin \overline{V}$ implies that $s\cdot x_0\notin \overline{V}\cdot x_0$. So, $r_{x_0}F'\in C_c^+({\cal G})$ and hence by taking a sequence  $(F_n)_{n\in \Bbb N}\subseteq  C_b({\cal X})$ for which $(r_{x_0}F_n)_{n\in \Bbb N}$ tends strictly to $g$, we can write $\|(r_{x_0}F_n)(r_{x_0}F')-g(r_{x_0}F')\|_\infty\longrightarrow 0$. Consequently, we have
\begin{eqnarray*}
\int_{\cal G} g(s)\,d\mu(s) &=& \int_{\cal G} g(s)(r_{x_0}F')(s)\,d\mu(s)\\
&=&  \langle\mu ,(r_{x_0}F')g\rangle\\
&=& \lim_{n\rightarrow +\infty} \langle\mu , (r_{x_0}F_n)(r_{x_0}F')\rangle \\
&=& \lim_{n\rightarrow +\infty} \int_{\cal G} F_n(s\cdot x_0)F'(s\cdot x_0)\,d\mu(s)\\
&=& \lim_{n\rightarrow +\infty} \int_{\cal G} (F_nF')(s\cdot x_0)\,d\mu(s).
\end{eqnarray*}
As $F_nF'\in C_b({\cal X})$, by using (\ref{rm2}), we deduce
$\int_{\cal G} g(s)\,d\mu(s)=0$. Now, since $g$ is arbitrary in $C_c({\cal G})$, we can deduce that $\mu=0$ and this completes the proof of this part.

%%%%%%%%%%%%%%%%%%%%%%%%%%%%%%%%%%%%%%%%%%%%%%%%%%%%%%%%%%%%%%%%%%%%%%%%%%%%%%%%%%%%%%%%

{\bf (b)} Let $\mu$ be a measure in $M({\cal G})$ such that
$\mu\star \sigma=0$ for all $\sigma\in M({\cal X})$.
Recall from \cite[Lemma 2.1. and Lemma 2.2.]{ghweight}
that the algebra $L^1({\cal G})$ possesses a
bounded approximate identity $(\varphi_\alpha)$
such that each $\varphi_\alpha$ is a continuous function
with compact support and $\lambda_{\varphi_\alpha}\longrightarrow\delta_e$
in the weak$^*$ topology of $M({\cal G})$.
As $M_a({\cal G})$ is an
ideal of $M({\cal G})$, we have $\lambda_{\varphi_\alpha}\ast\mu\in
M_a({\cal G})$ and also $$(\lambda_{\varphi_\alpha}\ast\mu)\star
\sigma=\lambda_{\varphi_\alpha}\star (\mu\star \sigma)=0,$$ for all
$\sigma\in M({\cal X})$ and all $\alpha$. Accordingly, from the assumption, we get
$\lambda_{\varphi_\alpha}\ast\mu=0$ for all $\alpha$. It follows that $\mu=0$.$\e$\\

%%%%%%%%%%%%%%%%%%%%%%%%%%%%%%%%%%%%%%%%%%%%%%%%%%%%%%%%%%%%%%%%%%%%%%%%%%%%%%%%%%%%%%%%%%%%%%%

When $\cal X$ is a $\cal G$-space, by the use of the action
of $\cal G$ on $\cal X$, we introduce a certain subspace of
$C_b({\cal X})$ which is the main object of study of this work. In
details, following Lau and Chu \cite{chulau} for a given $s$ in
$\cal G$, we define the left translation of a function $F\in
C_b({\cal X})$, by an element $s\in G$, by $l_sF(x)=F(s\cdot x)$ for all $x\in {\cal X}$
and consider
$$LUC({\cal X},{\cal G})=\Big{\{}F\in C_b({\cal X}):~ s\in {\cal G} \mapsto l_sF\in
(C_b({\cal X}),\|\cdot\|_\infty)  {\hbox{~is~
continuous}}\Big{\}}.$$ Then $LUC({\cal X},{\cal G})$ is a closed
subspace of $C_b({\cal X})$ which is invariant under left
translation. We always let $\cal G$ act on itself by left
translation and therefore $LUC({\cal G})=LUC({\cal G},{\cal G})$.
It is well-known that $C_0({\cal G})\subseteq LUC({\cal G})$ (see
\cite[Proposition 2.6]{farmonic}), by the same argument one can
prove the following result, which guarantees that $LUC({\cal
X},{\cal G})$ contains enough elements to separate the points of
$\cal X$.

\begin{lemma}\label{luc0}
Let $\cal G$ be a locally compact group and $\cal X$ be a $\cal
G$-space. Then the Banach space $C_0({\cal X})$ is contained in
$LUC({\cal X},{\cal G})$.
\end{lemma}

%%%%%%%%%%%%%%%%%%%%%%%%%%%%%%%%%%%%%%%%%%%%%%%%%%%%%%%%%%%%%%%%%%%%%%%%%%%%%%%%%%%%%%

Moreover, if $LUC({\cal X},{\cal G})^*$ denotes the first dual
space of the Banach space $LUC({\cal X},{\cal G})$, $M$ is an
arbitrary element of $LUC({\cal X},{\cal G})^*$ and $F\in
LUC({\cal X},{\cal G})$, we define the function $MF:{\cal
G}\rightarrow \Bbb{C}$ by
\begin{eqnarray*}\label{maction}
MF(s)=\langle M,l_sF\rangle ,
\end{eqnarray*}
for all $s\in {\cal G}$, which belongs to $LUC({\cal G})$; This is
because of,
\begin{eqnarray}\label{211228}
{_s}(MF)(t)=(MF)(st)
=\langle M,l_{st}F\rangle
=\langle M,l_t(l_sF)\rangle
=(M(l_sF))(t),
\end{eqnarray}
for all $s, t\in {\cal G}$. Also we have the following lemma which
plays a key role in this paper.

%%%%%%%%%%%%%%%%%%%%%%%%%%%%%%%%%%%%%%%%%%%%%%%%%%%%%%%%%%%%%%%%%%%%%%%%%%%%%%%%%%%%%%%%%%%%%%%

\begin{lemma}\label{1}
Let $\cal G$ be a locally compact group and $\cal X$ be a $\cal
G$-space. Then the mapping
\begin{eqnarray*}
\left\{\begin{array}{clcr}
LUC({\cal X},{\cal G})^*\times LUC({\cal X},{\cal G})\rightarrow LUC({\cal G})\\
(M,F)\mapsto MF
\end{array}
\right.
\end{eqnarray*}
is a bounded bilinear map with
$\|MF\|_\infty\leq\|M\|\;\|F\|_\infty$.
\end{lemma}

%%%%%%%%%%%%%%%%%%%%%%%%%%%%%%%%%%%%%%%%%%%%%%%%%%%%%%%%%%%%%%%%%%%%%%%%%%%%%%%%%%%%%%%%%%%%%%%

Lemma~\ref{1} paves the way for defining a bounded bilinear
map as follows
\begin{eqnarray*}
\left\{\begin{array}{clcr}
LUC({\cal G})^*\times LUC({\cal X},{\cal G})^*\rightarrow LUC({\cal X},{\cal G})^*\\
(m,M)\mapsto m\cdot M
\end{array}
\right.
\end{eqnarray*}
with $\|m\cdot M\|\leq\|m\|\;\|M\|$, where
\begin{eqnarray}\label{eq1}
\langle m\cdot M,F\rangle=\langle m,MF\rangle,
\end{eqnarray}
for all $F\in LUC({\cal X},{\cal G})$.

%%%%%%%%%%%%%%%%%%%%%%%%%%%%%%%%%%%%%%%%%%%%%%%%%%%%%%%%%%%%%%%%%%%%%%%%%%%%%%%%%%%%%%%%%%%%%%%

\begin{proposition}\label{2}
Let $\cal G$ be a locally compact group and $\cal X$ be a $\cal
G$-space. Under the mapping $(m,M)\mapsto m\cdot M$, defined by
{\rm(\ref{eq1}\rm)}, $LUC({\cal X},{\cal G})^*$ becomes a Banach
left $LUC({\cal G})^*$-module with $\|m\cdot M\|\leq\|m\|\;\|M\|$
and $\delta_e\cdot M=M$.
\end{proposition}
{\noindent Proof.} It is enough to prove that $(m\odot n)\cdot
M=m\cdot(n\cdot M)$, where $m, n\in LUC({\cal G})^*$ and $M\in
LUC({\cal X},{\cal G})^*$. To this end, suppose that $F\in
LUC({\cal X},{\cal G})$. Then, by (\ref{arensluc}) and (\ref{211228}), we can see that for all $s\in {\cal G}$
\begin{eqnarray*}
(n\diamond (MF))(s)=\langle n,{_s(MF)}\rangle =\langle
n,M(l_sF)\rangle =\langle n\cdot M,l_sF\rangle =((n\cdot M)F)(s),
\end{eqnarray*}
which shows that $n\diamond (MF)= (n\cdot M)F$. Therefore,
\begin{eqnarray*}
\langle(m\odot n)\cdot M,F\rangle &=&\langle m\odot n,MF\rangle\\
&=&\langle m,n\diamond (MF)\rangle\\
&=&\langle m,(n\cdot M)F\rangle\\
&=&\langle m\cdot(n\cdot M),F\rangle,
\end{eqnarray*}
and this completes the proof.$\e$\\

%%%%%%%%%%%%%%%%%%%%%%%%%%%%%%%%%%%%%%%%%%%%%%%%%%%%%%%%%%%%%%%%%%%%%%%%%%%%%%%%%%

If now, for a given $\sigma\in M({\cal X})$, we define a linear
functional on $LUC({\cal X},{\cal G})$, denoted again by $\sigma$,
which assigns to each $F\in LUC({\cal X},{\cal G})$ the value
$\int_{\cal X} F(x)\,d\sigma(x)$. Then, $M({\cal X})$ may be
regarded as a subspace of $LUC({\cal X},{\cal G})^*$. Moreover, it
is not hard to check that the inclusion $LUC({\cal G})^*\cdot
M({\cal X})\subseteq M({\cal X})$ can fail even if ${\cal X}={\cal
G}$, where
$$LUC({\cal G})^*\cdot M({\cal X})={\Big\{}m\cdot\sigma:~m\in LUC({\cal G})^*,
~\sigma\in M({\cal X}){\Big\}};$$
 In other word, the Banach space $M({\cal X})$ is not in
general an $LUC({\cal G})^*$-submodule of $LUC({\cal X},{\cal
G})^*$. On the other hand, if $\cal X$ is compact, then $M({\cal
X})$ is an $LUC({\cal G})^*$-submodule of $LUC({\cal X},{\cal
G})^*$; This is because of, in this case $M({\cal X})=LUC({\cal
X},{\cal G})^*$. We do not know if the converse of this fact is
valid in general; here, we prove the converse under an extra assumption. The
notation ${\cal CLS}({\cal X},{\cal G})$ in this proposition
and in the sequel denotes the norm closure of the linear span of
the set
$$LUC({\cal X},{\cal G})^*LUC({\cal X},{\cal G}):=
{\Big\{}MF:~ M\in LUC({\cal X},{\cal G})^*,~ F\in LUC({\cal
X},{\cal G}){\Big\}},$$ with respect to the norm topology of
$LUC({\cal G})$. Also, we say that the action of $LUC({\cal G})^*$
on $LUC({\cal X},{\cal G})^*$ is faithful if $m\in LUC({\cal
G})^*$ is so that $m\cdot M=0$ for all $M\in LUC({\cal X},{\cal
G})^*$, then $m=0$.

\begin{proposition}\label{400}
Let $\cal G$ be a locally compact group and $\cal X$ be a $\cal
G$-space. Then the following statements hold.
\newcounter{j12}
\begin{list}%
{\bf(\alph{j12})}{\usecounter{j12}}
\item The action of
$LUC({\cal G})^*$ on $LUC({\cal X},{\cal G})^*$ is faithful if and only if ${\cal CLS}({\cal
X},{\cal G})=LUC({\cal G})$.

\item If the action of $LUC({\cal G})^*$ on $LUC({\cal X},{\cal
G})^*$ is faithful, then $\cal G$ acts effectively on $\cal X$.

\item If $\cal X$ is compact, then $M({\cal X})$ is an $LUC({\cal
G})^*$-submodule of $LUC({\cal X},{\cal G})^*$. The converse is
also true if $\cal X$ is a transitive $\cal G$-space.
\end{list}
\end{proposition}

{\noindent Proof.} {\bf(a)} The proof of this assertion uses only
the Hahn-Banach theorem and so the details are omitted.

{\bf (b)} If the action of $LUC({\cal G})^*$ on $LUC({\cal
X},{\cal G})^*$ is faithful and $s\in {\frak N}({\cal X},{\cal
G})$, then for all $M\in LUC({\cal X},{\cal G})^*$ and $F\in
LUC({\cal X},{\cal G})$ we have
\begin{eqnarray*}
\langle \delta_{s}\cdot M , F\rangle =\langle \delta_{s},
MF\rangle =
 MF(s)=\langle M,l_sF\rangle= \langle M , F\rangle.
\end{eqnarray*}
So, $\delta_{s}\cdot M=M$ for all $M\in LUC({\cal X},{\cal
G})^*$ and hence, $\delta_{s}=\delta_e$; that is, $s=e$.

{\bf (c)} We only need to prove the converse of this assertion, which is the
essential part of it. To this end, suppose that $M({\cal X})$ is a
submodule of $LUC({\cal X},{\cal G})^*$. Then for all $m\in
LUC({\cal G})^*$ and $\sigma\in M({\cal X})$, we have $m\cdot
\sigma\in M({\cal X})$. This implies that the norm of the linear
functional $m\cdot \sigma$ on $LUC({\cal X},{\cal G})$ is obtained by
the following equality
$$\|m\cdot \sigma\|=\sup{\Big\{}|\langle m\cdot \sigma, F_0 \rangle|:\, F_0\in C_0({\cal
X}),\,\|F_0\|_\infty=1{\Big\}}.$$ Therefore, $$|\langle m\cdot
\sigma, F \rangle|\leq \sup {\Big\{}|\langle m\cdot \sigma, F_0
\rangle|:\, F_0\in C_0({\cal X}),\,\|F_0\|_\infty=1{\Big\}},$$ for
all $F\in LUC({\cal X},{\cal G})$ with $\|F\|_\infty=1$. Hence, for
all $m\in LUC({\cal G})^*$, $\sigma\in M({\cal X})$, and $F\in
LUC({\cal X},{\cal G})$ with $\|F\|_\infty=1$ we have
\begin{eqnarray}\label{eq3}
|\langle m, \sigma F  \rangle|\leq \sup {\Big\{}|\langle m, \sigma
F_0 \rangle|:\, F_0\in C_0({\cal X}),\,\|F_0\|_\infty=1{\Big\}}.
\end{eqnarray}
From this, we deduce that, for all $\sigma\in M({\cal X})$ and
$F\in LUC({\cal X},{\cal G})$, $\sigma F$ is in the norm closure
of the linear span of the set
$$\sigma C_0({\cal X}):={\Big\{}\sigma F_0:~ F_0\in C_0({\cal
X}){\Big\}},$$ with respect to the norm topology of $LUC({\cal
G})$; since, otherwise, we can find a $m\in LUC({\cal G})^*$ so that
$\langle m, \sigma F\rangle=1$ and $\langle m, \sigma
F_0\rangle=0$ for all $F_0\in C_0({\cal X})$, which is a
contradiction to (\ref{eq3}). In particular, if $x$ and $F$ are
arbitrary elements of $\cal X$ and $LUC({\cal X},{\cal G})$,
respectively, then, there exists a sequence
$(F_n)_{n\in\Bbb{N}}\subseteq C_0({\cal X})$ so that $\|\delta_x
F_n-\delta_x F\|_\infty\longrightarrow 0$. Now from the
transitivity of the action of $\cal G$ on $\cal X$, we get
\begin{eqnarray*}
\|F_n-F\|_\infty&\leq&\sup_{s\in {\cal G}}
|F_n(s\cdot x)-F(s\cdot x)|\\&\leq&\|\delta_x F_n-\delta_x F\|_\infty\\
&=&\sup_{s\in {\cal G}}
|(\delta_x F_n)(s)-(\delta_x F)(s)|.
\end{eqnarray*}
Thus, $\|F_n-F\|_\infty\longrightarrow 0$, which implies that $F$
is in $C_0({\cal X})$. Therefore, $LUC({\cal X},{\cal G})\subseteq
C_0({\cal X})$. As $LUC({\cal X},{\cal G})$ contains the constant
functions on $\cal X$,
 we deduce that $\cal X$ is compact.$\e$\\

%%%%%%%%%%%%%%%%%%%%%%%%%%%%%%%%%%%%%%%%%%%%%%%%%%%%%%%%%%%%%%%%%%%%%%%%%%%%%%%%%%
%%%%%%%%%%%%%%%%%%%%%%%%%%%%%%%%%%%%%%%%%%%%%%%%%%%%%%%%%%%%%%%%%%%%%%%%%%%%%%%%%%

\section{ The weak$^*$ continuity of the left $LUC({\cal G})^*$-module action }

As we know, Lau \cite{lau} has shown that a left translation
$n\mapsto m\odot n$ is weak$^*$ continuous on $LUC({\cal G})^*$
for a fixed $m$ in $LUC({\cal G})^*$ if and only if $m$ is in
$M({\cal G})$; that is, ${\frak{Z}}({\cal G})$, the topological
centre of $LUC({\cal G})^*$, is $M({\cal G})$. It follows that, if
$\cal X$ is a $\cal G$-space, then for an arbitrary $m$ in
$LUC({\cal G})^*$ the weak$^*$ continuity of the map $M\mapsto
m\cdot M$ on $LUC({\cal X},{\cal G})^*$ can fail even if ${\cal
X}={\cal G}$. A problem which is of interest is that for which
element $m\in LUC({\cal G})$ the map $M\mapsto m\cdot M$ on
$LUC({\cal X},{\cal G})^*$ is weak$^*$ to weak$^*$ continuous?
Therefore, it seems valuable to define
\begin{eqnarray*}
{{\frak{Z}}({\cal X},{\cal G})}=\Big{\{}m\in LUC({\cal
G})^*:\,M\mapsto m\cdot M~ {\hbox{is weak}}^*~{\hbox{to weak}}^*
~{\hbox{continuous on}}~LUC({\cal X},{\cal G})^*\Big{\}},
\end{eqnarray*}
the topological centre of the module action induced by $LUC({\cal
G})^*$ on $LUC({\cal X},{\cal G})^*$. In the special case that we
let $\cal G$ act on itself by left multiplication, the set ${\frak{Z}}({\cal G},{\cal G})$ coincides with
${\frak{Z}}({\cal G})$, which defined in Section 1.

%%%%%%%%%%%%%%%%%%%%%%%%%%%%%%%%%%%%%%%%%%%%%%%%%%%%%%%%%%%%%%%%%%%%%%%%%%%%%%%%%%%%%%%%%%%%%%%

This section studies the subspace ${{\frak{Z}}({\cal X},{\cal
G})}$ of $LUC({\cal G})^*$ in the case where $\cal X$ is a $\cal
G$-space and, in particular, the question when the subspace
${{\frak{Z}}({\cal X},{\cal G})}$ is $M({\cal G})$ or $LUC({\cal G})^*$. Before
proceeding further in this section, we should note that if $\cal
X$ is a $\cal G$-space, then $LUC({\cal X},{\cal G})^*$ is a left
Banach $G$-module. In fact, it is suffices to define the left action
of $\cal G$ on $LUC({\cal X},{\cal G})^*$ by $(s,M)\mapsto
\delta_s\cdot M$.

\begin{lemma}\label{30}
Let $\cal G$ be a locally compact group and $\cal X$ be a $\cal
G$-space. Then the action of $\cal G$ on the unit ball of
$LUC({\cal X},{\cal G})^*$ with the weak$^*$ topology is jointly continuous.
\end{lemma}
{\noindent Proof.} For the proof, it is suffices to note that if
$(s_\alpha)$ is a net converging to $s$ in $\cal G$, then the net
$(\delta_{s_\alpha})$ tends to $\delta_s$ with respect to the
weak$^*$ topology of $LUC({\cal G})^*$.$\e$\\

%%%%%%%%%%%%%%%%%%%%%%%%%%%%%%%%%%%%%%%%%%%%%%%%%%%%%%%%%%%%%%%%%%%%%%%%%%%%%%%%%%%%%%%%%%%%%%%

Let $\cal X$ be a $\cal G$-space. Given $F\in LUC({\cal X},{\cal
G})$ and $x\in{\cal X}$, we define $r_xF$ on $\cal G$ by
$$(r_xF)(s)=F(s\cdot x),\quad\quad\quad (s\in {\cal G}),$$ then a
routine computation shows that
$$\|{_{s_\alpha}}(r_xF)-{_{s}}(r_xF)\|_\infty\leq
\|l_{s_\alpha}F-l_{s}F\|_\infty,$$ where $(s_\alpha)$ is a net in
${\cal G}$ which tends to $s$, and therefore $r_xF$ is a function
in $LUC({\cal G})$. Hence, if $m$ is an arbitrary element of
$LUC({\cal G})^*$, then we can define a complex-valued function
$Fm$ on $\cal X$ by
$$Fm(x)=\langle m,r_xF\rangle\quad\quad (x\in {\cal X}).$$
Observe that $Fm$ is a bounded function on $\cal X$ with
$\|Fm\|_{\infty}\leq \|m\|\, \|F\|_{\infty}$. The following
lemma shows that the continuity of the function $Fm$ can fail
even if ${\cal X}={\cal G}$.

\begin{lemma}\label{9000}
Let $\cal G$ be a locally compact group and $\cal X$ be a $\cal
G$-space. Then, for all $m\in LUC({\cal G})^*$ and $F\in LUC({\cal X},{\cal G})$, the following assertions hold.
\newcounter{j14}
\begin{list}%
{\bf(\alph{j14})}{\usecounter{j14}}
\item For all $x\in {\cal X}$, $Fm(x)=\langle
m\cdot\delta_x,F\rangle$.
\item If $m\in
{\frak{Z}}({\cal X},{\cal G})$, then $Fm$ is in $C_b({\cal
X})$.
\item If $m\in {\frak{Z}}({\cal X},{\cal G})$, then
$\langle\delta_x,l_s(Fm)\rangle
=\langle(m\odot\delta_s)\cdot\delta_x,F\rangle$ for all $x\in
{\cal X}$ and $s\in {\cal G}$.
\end{list}
\end{lemma}
{\noindent Proof.} Let $m\in LUC({\cal G})^*$, $F\in LUC({\cal
X},{\cal G})$, $x\in {\cal X}$ and $s\in {\cal G}$ be given.

{\bf (a)} First, observe that
\begin{eqnarray}\label{210108}
(\delta_xF)(s)=\langle\delta_x, l_sF\rangle =l_sF(x)=F(s\cdot
x)=r_xF(s),
\end{eqnarray}
for all $s\in {\cal G}$. Hence, we have
\[\langle m\cdot\delta_x,F\rangle =\langle m,\delta_xF\rangle =
 \langle m,r_xF\rangle =Fm(x).\]

{\bf (b)} If $m\in {\frak{Z}}({\cal X},{\cal G})$ and
$(x_\alpha)\subseteq {\cal X}$ converges to $x\in{\cal X}$, then the net
$(\delta_{x_\alpha})$ tends to $\delta_x$ in the weak$^*$
topology of $LUC({\cal X},{\cal G})^*$. Therefore, since
$m\in{\frak{Z}}({\cal X},{\cal G})$, we have
\[(Fm)(x_\alpha)=\langle m\cdot\delta_{x_\alpha},F\rangle \longrightarrow
\langle m\cdot\delta_{x},F\rangle =(Fm)(x).\] This together with
the fact that $\|Fm\|_{\infty}\leq \|m\|\, \|F\|_{\infty}$,
implies that $Fm\in C_b({\cal X})$.

{\bf (c)} In this case, $r_xF$ is a function in $LUC({\cal G})$,
$Fm\in C_b({\cal X})$ and for each $t\in{\cal G}$ we have
\begin{eqnarray*}
\delta_s\diamond(r_xF)(t)=\langle\delta_s,{_t}(r_xF)\rangle =
{_t}(r_xF)(s)=F(t\cdot(s\cdot x))=(r_{s\cdot x}F)(t),
\end{eqnarray*}
From this, by (\ref{210108}), we deduce that
\begin{eqnarray*}
\langle\delta_x,l_s(Fm)\rangle=Fm(s\cdot x)=\langle m,r_{s\cdot x}F\rangle=
\langle m\odot\delta_s,r_xF\rangle
=\langle(m\odot\delta_s)\cdot\delta_x,F\rangle ,
\end{eqnarray*}
as desired.$\e$\\

%%%%%%%%%%%%%%%%%%%%%%%%%%%%%%%%%%%%%%%%%%%%%%%%%%%%%%%%%%%%%%%%%%%%%%%%%%%%%%%%%%%%%%

Now, suppose that $m$ is an element of $LUC({\cal G})^*$ such that,
for each $F\in LUC({\cal X},{\cal G})$, the function $Fm$ is in
$LUC({\cal X},{\cal G})$. Then every $M$ in $LUC({\cal X},{\cal
G})^*$ gives a linear functional $M\bullet m$ on $LUC({\cal
X},{\cal G})$ as follows
$$\langle M\bullet m, F\rangle =\langle M, Fm\rangle .$$
In addition, if $M\bullet m=m\cdot M$, for all $M\in LUC({\cal
X},{\cal G})^*$, then $m\in {\frak{Z}}({\cal X},{\cal G})$;
Indeed, if $(M_\alpha)\subseteq LUC({\cal X},{\cal G})^*$ tends to
an element $M$ with respect to the weak$^*$ topology, then we have
$$
\lim_\alpha\langle m\cdot M_\alpha, F\rangle= \lim_\alpha\langle
M_\alpha\bullet m, F\rangle=\langle M, Fm\rangle=\langle
M\bullet m, F\rangle=\langle m\cdot M, F\rangle,
$$
for all $F\in LUC({\cal X}, {\cal G})$. The following theorem
consider the converse of this fact whose proof is inspired by
\cite[Lemma 2.2]{lau}.

\begin{theorem}\label{5}
Let $\cal G$ be a locally compact group and $\cal X$ be a $\cal
G$-space. Then ${\frak{Z}}({\cal X},{\cal G})$ is precisely the
set of all $m\in LUC({\cal G})^*$ for which the following
conditions are satisfied
\newcounter{j13}
\begin{list}%
{\bf(\alph{j13})}{\usecounter{j13}} \item $Fm\in LUC({\cal
X},{\cal G})$ for all $F\in LUC({\cal X},{\cal G})$, \item
$M\bullet m=m\cdot M$ for all $M\in LUC({\cal X},{\cal G})^*$.
\end{list}
\end{theorem}

{\noindent Proof.} Suppose that $m\in {\frak{Z}}({\cal X},{\cal
G})$ and $F$ is an arbitrary element of
$LUC({\cal X},{\cal G})$. Then, the map $M\mapsto m\cdot M$ is weak$^*$ continuous on
$LUC({\cal X},{\cal G})^*$ and therefore it is weak$^*$ continuous
on bounded sets. Also, part {\bf(b)} of Lemma \ref{9000} implies that $Fm$ is a function in $C_b({\cal X})$.
Moreover, in light of Lemma \ref{30}, part {\bf (c)} of Lemma \ref{9000}, \cite[Lemma 2.5]{bami} and the weak$^*$ continuity of the map $M\mapsto m\cdot M$ on bounded set, we have
\begin{eqnarray}\label{210157}
\langle\Lambda,l_s(Fm)\rangle
=\langle(m\odot\delta_s)\cdot\Lambda,F\rangle
\end{eqnarray}
for all
$\Lambda\in C_b({\cal X})^*$ and $s\in {\cal G}$. Now, suppose
that $(s_\alpha)$ is a net in $\cal G$ converging to $s$ and for
each $\alpha$ the functional $\Lambda_\alpha\in C_b({\cal X})^*$
is chosen so that $\|\Lambda_\alpha\|=1$ and
$$\|l_{s_\alpha}(Fm)-l_s(Fm)\|_{\infty}={\Big|}
\langle\Lambda_\alpha,l_{s_\alpha}(Fm)-l_s(Fm)\rangle
{\Big|}.$$ Let also, $\Lambda$ be a weak$^*$-cluster point of
$(\Lambda_\alpha)$. By passing to a subnet, if necessary, we may assume
that $\Lambda_\alpha\longrightarrow\Lambda$ in the weak$^*$
topology of $LUC({\cal X},{\cal G})^*$.
Hence
\begin{eqnarray*}
\|l_{s_\alpha}(Fm)-l_s(Fm)\|_{\infty}&\leq&{\Big|}
\langle(m\odot\delta_{s_\alpha})\cdot\Lambda_\alpha-(m\odot\delta_{s})\cdot\Lambda,F\rangle {\Big|}\\&&
+{\Big|}\langle(m\odot\delta_{s})\cdot\Lambda-(m\odot\delta_{s})\cdot\Lambda_\alpha,F\rangle {\Big|}.
\end{eqnarray*}
Now, the weak$^*$ continuity of the map $M\mapsto m\cdot M$ on norm
bounded subsets of $LUC({\cal X},{\cal G})^*$ and Lemma \ref{30}
imply that $Fm$ is in $LUC({\cal X},{\cal G})$.

Finally, to prove the equality $M\bullet m=m\cdot M$ for all $M\in
LUC({\cal X},{\cal G})^*$. First note that $\delta_x\bullet
m=m\cdot \delta_x$, for all $x\in {\cal X}$. Hence, the equality
$M\bullet m=m\cdot M$ holds if $M$ is a convex combination of the
Dirac measures. We now invoke \cite[Lemma 2.5]{bami} to conclude
that $M\bullet m=m\cdot M$ holds for all $M\in LUC({\cal X},{\cal
G})^*$.$\e$\\

%%%%%%%%%%%%%%%%%%%%%%%%%%%%%%%%%%%%%%%%%%%%%%%%%%%%%%%%%%%%%%%%%%%%%%%%%%%%%%%%

It is clear that if $m\in {\frak Z}({\cal X},{\cal G})$, then
the map $M\mapsto m\cdot M$ is weak$^*$ continuous on all bounded
parts of $LUC({\cal X},{\cal G})^*$. As we have seen in the proof of Theorem \ref{5}, that
the map $M\mapsto m\cdot M$ is weak$^*$ continuous on all bounded
parts of $LUC({\cal X},{\cal G})^*$ plays a key role in its proof. Hence,
with an argument
similar to the proof of Theorem \ref{5} one
can prove the following result. The details are omitted.

\begin{corollary}\label{13-8-94}
Let $\cal G$ be a locally compact group and $\cal X$ be a $\cal
G$-space. Then $m\in {\frak Z}({\cal X},{\cal G})$ if and only if
the map $M\mapsto m\cdot M$ is weak$^*$ continuous on all bounded
parts of $LUC({\cal X},{\cal G})^*$.
\end{corollary}

%%%%%%%%%%%%%%%%%%%%%%%%%%%%%%%%%%%%%%%%%%%%%%%%%%%%%%%%%%%%%%%%%%%%%%%%%%%%%%%%%%%%%%%

The following two propositions paves the way for obtaining the wide class of
$\cal G$-space $\cal X$ for which $$M({\cal G})\subsetneqq{\frak Z}({\cal X},{\cal G})=LUC({\cal
G})^*.$$ First let us recall that, if $\cal H$ is a subgroup of $\cal G$, then the {\it
index} of $\cal H$ in $\cal G$, denoted by $[{\cal G}:{\cal H}]$,
is the number of left cosets of $\cal H$ in $\cal G$.

\begin{proposition}\label{0202}
Let $\cal G$ be a locally compact group and $\cal X$ be a
$\cal G$-space. If the index of ${\frak N}({\cal X},{\cal G})$ in
$\cal G$ is finite, then ${\frak Z}({\cal X},{\cal G})=LUC({\cal
G})^*$.
\end{proposition}
{\noindent Proof.} Let ${\frak N}:={\frak N}({\cal X},{\cal G})$
and $m$ be an arbitrary element of $LUC({\cal G})^*$. To prove
that $m\in {\frak Z}({\cal X},{\cal G})$, by Corollary
\ref{13-8-94} above, it will be enough to prove that the map
$M\mapsto m\cdot M$ is weak$^*$ continuous on bounded parts of
$LUC({\cal X},{\cal G})^*$. To this end, suppose that $(M_\alpha)$
is a bounded net in $LUC({\cal X},{\cal G})^*$ which tends to
$M\in LUC({\cal X},{\cal G})^*$ with respect to the weak$^*$
topology. Then, since the index of ${\frak
N}$ in $\cal G$ is finite, we can obtain a partition
${\Big\{}s_1{\frak N},s_2{\frak N},\cdots,s_k{\frak N}{\Big\}}$
for $\cal G$ and therefore for each $s\in{\cal G}$, there exists
$t\in{\frak N}$ and a unique $1\leq i\leq k$ for which $s=s_it$.
It follows that for all $\alpha$, $s\in{\cal G}$ and $F\in LUC({\cal X},{\cal
G})$
\begin{eqnarray*}
M_\alpha F(s)=\sum_{i=1}^kM_\alpha F(s_i)\chi_{s_i{\frak
N}}(s),\quad {\hbox{and}}\quad
MF(s)=\sum_{i=1}^kMF(s_i)\chi_{s_i{\frak N}}(s),
\end{eqnarray*}
 where
$\chi_{s_i{\frak N}}$ denotes the characteristic function of
${s_i{\frak N}}$ on $\cal G$.
On the other
hand, since, for each $F\in LUC({\cal X},{\cal G})$, $M_\alpha
F\longrightarrow MF$ pointwisely, we may choose $\alpha_0$ such
that $|M_\alpha F(s_i)-MF(s_i)|<\varepsilon$ whenever
$i=1,\cdots,k$ and $\alpha\succeq\alpha_0$. Therefore, for each
$\alpha\succeq\alpha_0$, we have
$$\|M_\alpha F-MF\|_\infty={\Big\|}\sum_{i=1}^kM_\alpha F(s_i)\chi_{s_i{\frak
N}}-\sum_{i=1}^kMF(s_i)\chi_{s_i{\frak N}}{\Big\|}_\infty<\varepsilon.$$
This implies that
$$|\big<m\cdot M_\alpha,F\big>-\big<m\cdot M,F\big>|\leq\|m\|\;\|
M_\alpha F-MF\|_\infty\longrightarrow 0.$$ Therefore, $m\in {\frak
Z}({\cal X},{\cal G})$.$\e$\\

%%%%%%%%%%%%%%%%%%%%%%%%%%%%%%%%%%%%%%%%%%%%%%%%%%%%%%%%%%%%%%%%%%%%%%%%%%%%%%%%%%%%%%%%%%

The following proposition, gives some of the main properties of
the set ${\frak{Z}}({\cal X},{\cal G})$.

\begin{proposition}\label{8}
Let $\cal G$ be a locally compact group and $\cal X$ be a $\cal
G$-space. Then the following assertions hold.
\newcounter{j2}
\begin{list}%
{\bf(\alph{j2})}{\usecounter{j2}} \item ${\frak{Z}}({\cal X},{\cal
G})$ is a subalgebra of $LUC({\cal G})^*$.

\item ${\frak{Z}}({\cal X},{\cal G})$ is closed with respect to
the norm topology of $LUC({\cal G})^*$.

\item $M({\cal G})$ is contained in ${\frak{Z}}({\cal X},{\cal
G})$.
\end{list}
\end{proposition}
{\noindent Proof}. In order to prove {\bf(a)}, assume that $m$ and
$n$ are in ${\frak{Z}}({\cal X},{\cal G})$. Suppose also that $F$
is an arbitrary element of $LUC({\cal X},{\cal G})$. Then, by Theorem \ref{5},
the functions $F':=Fm$ and $F'n$ are in $LUC({\cal X},{\cal G})$.
On the other hand, for each $x\in {\cal X}$, we have
\begin{eqnarray*}
(F'n)(x)&=&\langle n\cdot\delta_x,Fm\rangle\\
&=&\langle(n\cdot\delta_x)\bullet m,F\rangle\quad\quad\quad ({\hbox{since}}~
Fm\in LUC({\cal X},{\cal
G}))\\
&=&\langle m\cdot(n\cdot\delta_x),F\rangle~\quad\quad\quad ({\hbox{since}}~
m\in{\frak{Z}}({\cal X},{\cal
G}))\\
&=&\langle (m\odot n)\cdot\delta_x,F\rangle\\
&=&(F(m\odot n))(x).
\end{eqnarray*}
Hence, $F(m\odot n)$ is in $LUC({\cal X},{\cal G})$. Moreover,
for each $M\in LUC({\cal X},{\cal G})^*$, we have
\begin{eqnarray*}
\langle(m\odot n)\cdot M,F\rangle&=&\langle m\cdot(n\cdot M),F\rangle\\
&=&\langle (n\cdot M)\bullet
m,F\rangle\\
&=&\langle M\bullet n,Fm\rangle\\ &=&\langle M,
F'n\rangle\\&=&\langle M\bullet(m\odot n), F\rangle.
\end{eqnarray*}
Therefore, Theorem \ref{5} implies that
$m\odot n$ is in ${\frak{Z}}({\cal X},{\cal G})$.

Now, for the proof of the assertion {\bf(b)}, let $(m_k)$ be a
sequence in ${\frak{Z}}({\cal X},{\cal G})$ which tends to $m\in
LUC({\cal G})^*$ with respect to the norm topology. Also, suppose that
$F\in LUC({\cal X},{\cal G})$ and $M\in LUC({\cal X},{\cal G})^*$ are given.
From Lemma
\ref{9000}, we observe that
\begin{eqnarray*}
\|Fm_k-Fm\|_\infty&=&\sup_{x\in{\cal X}}|Fm_k(x)-Fm(x)|\\
&=&\sup_{x\in{\cal X}}|\langle(m_k-m)\cdot\delta_x,F\rangle |\\
&\leq&\|m_k-m\|\|F\|_\infty.
\end{eqnarray*}
This together with the fact that
 $m_k\in{\frak{Z}}({\cal X},{\cal G})$ implies that $Fm$ is an element of
$LUC({\cal X},{\cal G})$. On the other hand,
\begin{eqnarray*}
\langle m\cdot M,F\rangle
=\lim_k\langle M\bullet m_k, F\rangle
=\lim_k\langle M, Fm_k\rangle
=\langle M, Fm\rangle
=\langle M\bullet m, F\rangle.
\end{eqnarray*}
We now invoke Theorem \ref{5} to conclude that
$m\in {\frak{Z}}({\cal X},{\cal G})$.

%%%%%%%%%%%%%%%%%%%%%%%%%%%%%%%%%%%%%%%%%%%%%%%%%%%%%%%%%%%%%%%%%%%%%%%%%%%%%%%%%%%%%%%%%%%%%%%

Finally, for the proof of the last assertion, let $\mu$ be a
measure in $M({\cal G})$. Since the
measures in $M({\cal G})$ with compact supports are norm dense in
$M({\cal G})$, without loss of generality, we may
assume that $\mu$ has compact support. Let also $(M_\alpha)$ be a bounded net in
$LUC({\cal X},{\cal G})^*$ such that $M_\alpha\rightarrow M$ with
respect to the weak$^*$ topology, for some $M\in LUC({\cal
X},{\cal G})^*$. Choose $K>0$ such that $\|M_\alpha\|,\,\|M\|\leq
K$ for all $\alpha$. For each $F\in LUC({\cal X},{\cal G})$, $\alpha$ and
$s, t\in {\cal G}$ we have
\[|M_\alpha F(s)-M_\alpha F(t)|=|\langle M_\alpha,l_sF\rangle
-\langle M_\alpha ,l_tF\rangle |\leq K\|l_sF-l_tF\|_\infty.\] Hence, the
family of functions $M_\alpha F$ is equicontinuous. Therefore,
since $M_\alpha
F\longrightarrow MF$ pointwisely, the net
$M_\alpha F$  convergence uniformly to $MF$ on every compact
subset of $\cal G$. Thus,
\[|\langle\mu\cdot M_\alpha,F\rangle -\langle\mu\cdot M,F\rangle|
={\Big|}\int_{\cal G} M_\alpha F\;d\mu-\int_{\cal G}
MF\;d\mu{\Big|}\longrightarrow 0.\] It follows that the map
$M\mapsto\mu\cdot M$ is weak$^*$ continuous on all bounded parts
of $LUC({\cal X},{\cal G})^*$. We now invoke Corollary
\ref{13-8-94} to conclude that
$\mu\in {\frak{Z}}({\cal X},{\cal G})$.$\e$\\

%%%%%%%%%%%%%%%%%%%%%%%%%%%%%%%%%%%%%%%%%%%%%%%%%%%%%%%%%%%%%%%%%%%%%%%%%%%%%%%%%%%%%

The following result gives some necessary
 and sufficient conditions
for the validity of the equality ${\frak{Z}}({\cal X},{\cal
G})=M({\cal G})$.

\begin{theorem}\label{13}
Let $\cal G$ be a locally compact non-compact group and let $\cal
X$ be a $\cal G$-space. Then the following assertions are
equivalent.
\newcounter{j19}
\begin{list}%
{\bf(\alph{j19})}{\usecounter{j19}} \item ${\frak{Z}}({\cal
X},{\cal G})=M({\cal G})$. \item ${\cal CLS}({\cal X},{\cal
G})=LUC({\cal G})$. \item The action of $LUC({\cal G})^*$ on
$LUC({\cal X},{\cal G})^*$ is faithful.
\end{list}
\end{theorem}

{\noindent Proof.} By Proposition \ref{400}, we only need to prove
{\bf (a)}$\Leftrightarrow${\bf(b)}. In order
to prove that {\bf(a)} implies {\bf(b)},
suppose that ${\frak{Z}}({\cal X},{\cal G})=M({\cal G})$. If we
had an element $f$ in $LUC({\cal G})$ not belonging to the
subspace ${\cal CLS}({\cal X},{\cal G})$, then, by the Hahn-Banach
Theorem, we would have a functional $m\in LUC({\cal G})^*$
vanishing on the subspace ${\cal CLS}({\cal X},{\cal G})$ and such
that $\big<m,f\big>\neq 0$. From this, we deduce that $m\cdot M=0$
for all $M\in LUC({\cal X},{\cal G})^*$. It follows that $m\in
{\frak{Z}}({\cal X},{\cal G})=M({\cal G})$. Now, let
$$C_0({\cal G})^\perp={\Big\{}n\in LUC({\cal G})^*:~\big<n,f\big>=0~~{\hbox{for all}}~f\in C_0({\cal G}){\Big\}},$$
and pick any $0\neq m'\in C_0({\cal G})^\perp$
that is right cancellable in $LUC({\cal G})^*$ (such points exists
by \cite[Theorem 4]{filrc}). Note that for $m'$ and each $M\in
LUC({\cal X},{\cal G})^*$, we have
$$\big<(m\odot m')\cdot M,F\big>=\big<m\cdot(m'\cdot M),F\big>=0$$
for all $F\in LUC({\cal X},{\cal G})$. Therefore, $(m\odot
m')\cdot M=0$ for all $M\in LUC({\cal X},{\cal G})^*$ and this
implies that $m\odot m'\in {\frak{Z}}({\cal X},{\cal G})=M({\cal
G})$. From this, by \cite[Lemma 1.1]{glausert}, we deduce that
$$m\odot m'\in M({\cal G})\cap C_0({\cal G})^\perp=\{0\}.$$ This together
with the fact that $0\neq m'$ is right cancellable, implies that $m=0$ which is impossible.

To prove implication {\bf(b)}$\Rightarrow${\bf(a)}, suppose that
${\cal CLS}({\cal X},{\cal G})=LUC({\cal G})$. By part {\bf (c)}
of Proposition \ref{8}, the inclusion $M({\cal
G})\subseteq{\frak{Z}}({\cal X},{\cal G})$ being always true, it
will be enough to prove the reverse inclusion. To this end,
suppose that $m\in {\frak{Z}}({\cal X},{\cal G})$ and $(n_\alpha)$
is a bounded net in $LUC({\cal G})^*$ such that
$n_\alpha\longrightarrow n$ in $LUC({\cal G})^*$ with respect to
the weak$^*$-topology. Hence, for every $M\in LUC({\cal X},{\cal
G})^*$ and $F\in LUC({\cal X},{\cal G})$, we have
\begin{eqnarray*}
\langle m\odot n_\alpha,MF\rangle &=&\langle(m\odot n_\alpha)\cdot
M,F\rangle\\ &=&\langle m\cdot(n_\alpha\cdot M),F\rangle
\\&\longrightarrow &\langle m\cdot(n\cdot M),F\rangle\\ &=&\langle m\odot
n,MF\rangle ;
\end{eqnarray*}
This is because of, the net $(n_\alpha\cdot M)$ is a net in $LUC({\cal X},{\cal G})^*$
which converges to $nM$ in the weak$^*$ topology of $LUC({\cal X},{\cal G})^*$.
We now invoke the equality ${\cal CLS}({\cal X},{\cal
G})=LUC({\cal G})$ to conclude that
$$\langle m\odot n_\alpha,f\rangle \longrightarrow\langle m\odot n,f\rangle ,$$
for all $f\in LUC({\cal G})$. From this, by \cite[Theorem
4.1 and Lemma 2.2]{lau}, we deduce that $m$ belongs to ${\frak{Z}}({\cal
G})=M({\cal G})$. Therefore, ${\frak{Z}}({\cal X},{\cal
G})\subseteq M({\cal G})$. Hence, we have the equality
${\frak{Z}}({\cal X},{\cal G})=M({\cal G})$.$\e$\\

%%%%%%%%%%%%%%%%%%%%%%%%%%%%%%%%%%%%%%%%%%%%%%%%%%%%%%%%%%%%%%%%%%%%%%%%%%%%%%%%%%%%

We conclude this section with the following result which is a
consequence of Theorem \ref{210422}, part (b) of Proposition \ref{400}
and  Theorem~\ref{13}.

\begin{corollary}\label{500}
Let $\cal G$ be a locally compact non-compact group and $\cal X$
be a $\cal G$-space. If $\cal G$ does not act effectively on $\cal
X$, then $M({\cal G})$ is properly contained in ${\frak{Z}}({\cal
X},{\cal G})$.
\end{corollary}

%%%%%%%%%%%%%%%%%%%%%%%%%%%%%%%%%%%%%%%%%%%%%%%%%%%%%%%%%%%%%%%%%%%%%%%%%%%%%%%%%%%%%%
%%%%%%%%%%%%%%%%%%%%%%%%%%%%%%%%%%%%%%%%%%%%%%%%%%%%%%%%%%%%%%%%%%%%%%%%%%%%%%%%%%%%%%

\section{\normalsize\bf Examples}\label{s3}

This section is devoted to examples with two different purposes.
First, illustrating the results presented
in this paper for certain $\cal G$-spaces $\cal X$. Second, to give some example for which

$\centerdot$ ${\frak Z}({\cal X},{\cal G})=M({\cal G})$;

$\centerdot$ $M({\cal G})\subsetneqq{\frak Z}({\cal X},{\cal G})=LUC({\cal G})^*$;

$\centerdot$ ${\frak Z}({\cal X},{\cal G})$ is neither $M({\cal G})$ nor $LUC({\cal G})^*$.\\ To this end,
we commence this section with the following example which characterize the subalgebra
${\frak Z}({\cal X},{\cal G})$ of $LUC({\cal G})^*$ for certain $\cal G$-spaces $\cal X$.

%%%%%%%%%%%%%%%%%%%%%%%%%%%%%%%%%%%%%%%%%%%%%%%%%%%%%%%%%%%%%%%%%%%%%%%%%%%%%%%%%%%%%%%%%%%%%%%%%%%%%%%

\begin{example}\label{exmpl1112}
{\rm Let $\cal G$ be a locally compact group and $\cal X$ be a
$\cal G$-space.
\newcounter{q1021}
\begin{list}%
{\bf(\alph{q1021})}{\usecounter{q1021}} \item If ${\cal X}={\cal
G}$, then $M({\cal
G})={\frak{Z}}({\cal X},{\cal G})$.

\item If $\cal G$ is compact, then
\begin{eqnarray*}
M({\cal G})\subseteq {\frak{Z}}({\cal X},{\cal G})\subseteq
LUC(G)^*=M({\cal G})
\end{eqnarray*}
which shows that ${\frak{Z}}({\cal X},{\cal G})=M({\cal G})$.

\item If $\cal X$ is a finite discrete space, then
${\frak{Z}}({\cal X},{\cal G})=LUC({\cal G})^*$; Indeed, the
vector space $LUC({\cal X},{\cal G})$, as a subspace of $C_b({\cal
X})$, is a finite dimensional vector space, and so all locally
convex topologies on finite dimensional space $LUC({\cal X},{\cal G})^*$ are coincide. In
particular, for each $m\in LUC({\cal G})^*$ the linear map
$M\mapsto m\cdot M$ is continuous with respect to the norm
topology of $LUC({\cal X},{\cal G})^*$,
and also it is continuous with respect to the weak$^*$ topology. Therefore,
${\frak{Z}}({\cal X},{\cal G})=LUC({\cal G})^*$.
Specially, if $\cal G$ is non-compact, then $M({\cal G})\subsetneqq
{\frak Z}({\cal X},{\cal G})=LUC({\cal G})^*$.
\end{list}
}
\end{example}

%%%%%%%%%%%%%%%%%%%%%%%%%%%%%%%%%%%%%%%%%%%%%%%%%%%%%%%%%%%%%%%%%%%%%%%%%%%%%%%%%%%%%%%%%%%%%%%

It is worthwhile to mention that when $\cal G$ is a locally
compact group and $\cal H$ is a closed subgroup of $\cal G$, then
the space ${\cal G}/{\cal H}$ consisting of all left cosets of
$\cal H$ in $\cal G$ is a locally compact Hausdorff topological
space on which $\cal G$ acts from the left by
$${\cal G}\times{\cal G}/{\cal H}\rightarrow{\cal G}/{\cal H};~(s,t{\cal
H})\mapsto (st){\cal H}.$$ It has been shown
that if $\cal G$ is $\sigma$-compact, then every transitive $\cal
G$-space is homeomorphic to the quotient space ${\cal G}/{\cal H}$
for some closed subgroup $\cal H$, see
\cite[Subsection~2.6]{farmonic}.

%%%%%%%%%%%%%%%%%%%%%%%%%%%%%%%%%%%%%%%%%%%%%%%%%%%%%%%%%%%%%%%%%%%%%%%%%%%%%%%%%%%%%%%%%%%%%%%

In what follows, the notation ${\frak C}_{\cal G}({\cal H})$ is
used to denote the centralizer of the subgroup $\cal H$ in $\cal
G$; That is, $${\frak C}_{\cal G}({\cal H})={\Big\{}s\in {\cal G}:\,sh=hs,
\,{\hbox{for all}}~h\in {\cal H}{\Big\}}.$$

\begin{example}\label{exmpl111}
{\rm Let $\cal G$ be a locally compact, non-compact group and let
$\cal H$ and $\cal K$ be two non-trivial closed normal subgroups
of $\cal G$.
\newcounter{q1020}
\begin{list}%
{\bf(\alph{q1020})}{\usecounter{q1020}} \item If ${\cal X}={\cal
G}/{\cal H}$, then, by Corollary \ref{500}, $M({\cal
G})\varsubsetneq {\frak{Z}}({\cal X},{\cal G})$; This is because
of, in this case
\begin{eqnarray*}
{\frak N}({\cal X},{\cal G})={\Big \{}s\in {\cal G}:\,st{\cal H}=t{\cal
H},\,t\in {\cal G}{\Big\}}= {\cal H}.
\end{eqnarray*}
In particular, if the index of $\cal H$ in $\cal G$ is finite, then, by
Proposition \ref{0202}, we have ${\frak{Z}}({\cal X},{\cal
G})=LUC({\cal G})^*$.

\item If ${\cal X}={\cal K}$, then, obviously, $\cal G$ acts on
$\cal X$ by conjugation and makes ${\cal X}$ into a $\cal
G$-space; That is, $${\cal G}\times {\cal X}\rightarrow {\cal
X};~(s,x)\mapsto x^s:=sxs^{-1}.$$ In this case, ${\frak N}({\cal
X},{\cal G})={\frak C}_{\cal G}({\cal K})$. Hence, we can say if
either ${\frak C}_{\cal G}({\cal G})$ or ${\frak C}_{\cal G}({\cal
K})$ is non-trivial (for example, if $\cal K$ is abelian), then we
have $M({\cal G})\varsubsetneq {\frak{Z}}({\cal X},{\cal G})$.
\end{list}
}
\end{example}

%%%%%%%%%%%%%%%%%%%%%%%%%%%%%%%%%%%%%%%%%%%%%%%%%%%%%%%%%%%%%%%%%%%%%%%%%%%%%%%%%%%%%

Now, we consider the following special case of the previous example
which illustrates Proposition \ref{0202} and gives an example of $\cal G$ and $\cal X$
such that $M({\cal G})\subsetneqq
{\frak Z}({\cal X},{\cal G})=LUC({\cal G})^*$.

\begin{example}
{\rm Let $\Bbb Z$ be the discrete group of integer numbers and
let, for each $n\in {\Bbb Z}^+$, $n{\Bbb Z}$ be the subgroup generated
by $n$ consists of all integer multiples of $n$. If ${\cal X}_n={\Bbb Z}/{n{\Bbb Z}}$,
then Corollary \ref{500} together with the fact that
${\frak N}({\cal X}_n,{\Bbb Z})=n{\Bbb Z}$ implies that
$M({\Bbb Z})\varsubsetneq {\frak{Z}}({\cal X}_n,{\Bbb Z})$.
On the other hand, since the index of $n{\Bbb Z}$ in $\Bbb Z$ is $n$, by
Proposition \ref{0202}, we have ${\frak{Z}}({\cal X}_n,{\Bbb Z})=LUC({\Bbb Z})^*$.}
\end{example}

%%%%%%%%%%%%%%%%%%%%%%%%%%%%%%%%%%%%%%%%%%%%%%%%%%%%%%%%%%%%%%%%%%%%%%%%%%%%%%%%%%%%%%%%%

The following example illustrates Theorem \ref{13}.

\begin{example}
{\rm Suppose that $\cal Y$ is an arbitrary locally compact
Hausdorff space and $\cal G$ is a locally compact
group. Then  ${\cal X}={\cal G}\times {\cal Y}$, equipped with the
product topology and the action defined by
$${\cal G}\times{\cal X}\rightarrow{\cal X};~s\cdot(t,y)=(st,y),$$
is a $\cal G$-space for which ${\frak Z}({\cal X},{\cal G})=M({\cal G})$;
Indeed, by part (b) of Example \ref{exmpl1112}, we only need to show that the equality
${\frak Z}({\cal X},{\cal G})=M({\cal G})$
is valid when $\cal G$ is non-compact. To this end, first note that,
if $f\in LUC({\cal G})$ and $F\in C_b({\cal Y})$, then
the function $$(f\otimes F):{\cal X}\rightarrow{\Bbb C};~(t,y)\mapsto f(t)F(y),$$
is a function in $LUC({\cal X},{\cal G})$; This is because of,
$$l_s(f\otimes F)=({}_sf)\otimes F,$$ for all
$s\in {\cal G}$. In particular, if $1_{\cal Y}$ denotes the constant function 1 on $\cal Y$, then, for each
$f\in LUC({\cal G})$, the function $f\otimes 1_{\cal Y}$ is a function in $LUC({\cal X},{\cal G})$ and therefore we can consider $LUC({\cal G})$ as a closed subspace of $LUC({\cal X},{\cal G})$ via the inclusion mapping
$$\iota:LUC({\cal G})\rightarrow LUC({\cal X},{\cal G});~f\mapsto f\otimes 1_{\cal Y}.$$
Then $\iota^*:LUC({\cal X},{\cal G})^*\rightarrow LUC({\cal G})^*$ is the restriction mapping and hence norm decreasing and onto by Hahn-Banach Theorem.
Moreover, if $m \in LUC({\cal G})^*$ and $M\in
LUC({\cal X},{\cal G})^*$ is chosen so that $\iota^*(M)=m$, then
\begin{eqnarray*}
(m\diamond f)(s)&=&\langle m,{}_sf\rangle\\&=&\langle \iota^*(M),{}_sf\rangle\\&=&
\langle
M,({}_sf)\otimes1_{\cal Y}\rangle\\&=&\langle M,l_s(f\otimes1_{\cal Y})\rangle\\&=&(M(f\otimes1_{\cal Y}))(s),
\end{eqnarray*}
for all $f\in LUC({\cal G})$ and $s\in {\cal G}$.
Hence, we have
\begin{eqnarray*}
{\Big\{}m\diamond f:\,m\in LUC({\cal G})^*,\, f\in LUC({\cal
G}){\Big\}} \subseteq {\Big\{}MF:\,M\in LUC({\cal X},{\cal
G})^*,\,F\in LUC({\cal X},{\cal G}){\Big\}},
\end{eqnarray*}
which shows that ${\cal CLS}({\cal X},{\cal G})=LUC({\cal G})$. We now invoke Theorem \ref{13} to conclude that ${\frak{Z}}({\cal X},{\cal G})=M({\cal G})$.}
\end{example}

The proof of Example \ref{0925} below, which gives a $\cal G$-space $\cal X$
such that ${\frak N}({\cal X},{\cal G})$ is neither $M({\cal G})$ nor $LUC({\cal G})^*$, relies on the following example.

\begin{example}\label{rem1}
{\rm Let ${\cal G}_1$ and ${\cal G}_2$ be two locally compact groups. Then we can
consider $LUC({\cal G}_1)$ as a closed subspace of
$LUC({\cal G}_1\times {\cal G}_2)$ via the inclusion mapping
$$\iota_{{\cal G}_1}:LUC({\cal G}_1)\rightarrow LUC({\cal G}_1\times {\cal G}_2);~
f\mapsto f\otimes 1_{{\cal G}_2},$$
where $1_{{\cal G}_2}$ denotes the constant function 1 on ${{\cal G}_2}$.
Suppose also that
${\cal X}_1$ is a ${\cal G}_1$-space and ${\cal X}_2$ is
a ${\cal G}_2$-space. Then
${\cal X}_1\times {\cal X}_2$, with product topology and coordinatewise
left module action, is a ${\cal G}_1\times {\cal G}_2$-space for which
$$
\iota_{{\cal G}_1}^*{\Big(}{\frak{Z}}({\cal X}_1\times {\cal X}_2, {\cal G}_1\times
{\cal G}_2){\Big)}\subseteq{\frak{Z}}({\cal X}_1, {\cal G}_1).
$$
In particular, since the map $\iota_{{\cal G}_1}^*$ is onto, if $${\frak{Z}}({\cal X}_1\times {\cal X}_2, {\cal G}_1\times
{\cal G}_2)=LUC({\cal G}_1\times{\cal G}_2)^*,$$ then ${\frak{Z}}({\cal X}_1, {\cal G}_1)=LUC({\cal G}_1)^*$; Indeed,
if $m\in {\frak{Z}}({\cal X}_1\times {\cal X}_2, {\cal G}_1\times
{\cal G}_2)$, $F_1\in LUC({\cal X}_1,{\cal G}_1)$ and $1_{{\cal X}_2}$ denotes the constant function on ${\cal X}_2$ with value 1, then
$F_1\otimes 1_{{\cal X}_2}\in LUC({\cal X}_1\times {\cal X}_2, {\cal G}_1\times
{\cal G}_2)$ and we have
\begin{eqnarray*}
{\Big(}(F_1\otimes 1_{{\cal X}_2})m{\Big)}(x_1,x_2) &=& \langle m,r_{(x_1,x_2)}(F_1\otimes 1_{{\cal X}_2})\rangle\\
&=& \langle m,(r_{x_1}F_1)\otimes 1_{{\cal X}_2})\rangle\\
&=& \langle \iota_{{\cal G}_1}^*(m),r_{x_1}F_1\rangle\\
&=& {\Big(}F_1\iota_{{\cal G}_1}^*(m){\Big)}(x_1),
\end{eqnarray*}
for all $x_1\in {\cal X}_1$ and $x_2\in {\cal X}_2$ and therefore
$$(F_1\otimes 1_{{\cal X}_2})m=(F_1\iota_{{\cal G}_1}^*(m))\otimes 1_{{\cal X}_2}.$$
From this, by Theorem \ref{5}, we have
$$(F_1\iota_{{\cal G}_1}^*(m))\otimes 1_{{\cal X}_2}\in LUC({\cal X}_1\times {\cal X}_2, {\cal
G}_1\times {\cal G}_2).$$ It follows that $F_1\iota_{{\cal G}_1}^*(m)\in LUC({\cal X}_1,{\cal
G}_1)$ for all $F_1\in LUC({\cal X}_1,{\cal G}_1)$. Moreover, if
$$\iota_{{\cal X}_1}:LUC({\cal X}_1,{\cal G}_1)\rightarrow LUC({\cal X}_1\times {\cal X}_2,{\cal G}_1\times {\cal G}_2);~
F_1\mapsto F_1\otimes 1_{{\cal X}_2},$$
then it is not hard to check that
$$M(F_1\otimes 1_{{\cal X}_2})=(\iota_{{\cal X}_1}^*(M)F_1)\otimes 1_{{\cal G}_2},$$
for all $F_1\in LUC({\cal X}_1,{\cal G}_1)$, where $1_{{\cal G}_2}$ denotes
the constant function one on ${\cal G}_2$. It follows that
$$\iota_{{\cal G}_1}^*(m)\cdot \iota_{{\cal X}_1}^*(M)=\iota_{{\cal X}_1}^*(M)\bullet \iota_{{\cal G}_1}^*(m),$$
for all $M\in LUC({\cal X}_1\times {\cal X}_2,{\cal G}_1\times {\cal G}_2)^*$; In details, if $F_1\in LUC({\cal X}_1,{\cal G}_1)$, then, we can write
\begin{eqnarray*}
\langle \iota_{{\cal G}_1}^*(m)\cdot \iota_{{\cal X}_1}^*(M), F_1\rangle
&=& \langle \iota_{{\cal G}_1}^*(m), \iota_{{\cal X}_1}^*(M)F_1\rangle\\
&=& \langle m,(\iota_{{\cal X}_1}^*(M)F_1)\otimes 1_{{\cal G}_2}\rangle\\
&=& \langle m,M (F_1\otimes 1_{{\cal X}_2})\rangle\\
&=& \langle m\cdot M ,F_1\otimes 1_{{\cal X}_2}\rangle\\
&=& \langle M\bullet m ,F_1\otimes 1_{{\cal X}_2}\rangle~\quad\quad\quad{\Big(}{\hbox{since}}~ m\in {\frak{Z}}({\cal X}_1\times {\cal X}_2, {\cal G}_1\times
{\cal G}_2){\Big)}\\
&=& \langle M, (F_1\otimes 1_{{\cal X}_2})m\rangle\\
&=& \langle M, (F_1\iota_{{\cal G}_1}^*(m))\otimes 1_{{\cal X}_2}\rangle\\
&=& \langle \iota_{{\cal X}_1}^*(M), F_1\iota_{{\cal G}_1}^*(m)\rangle\\
&=& \langle \iota_{{\cal X}_1}^*(M)\bullet \iota_{{\cal G}_1}^*(m), F_1\rangle.\quad
\quad{\Big(}{\hbox{since}}~ F_1\iota_{{\cal G}_1}^*(m)\in LUC({\cal X}_1,{\cal
G}_1){\Big)}
\end{eqnarray*}
Hence, since the map $\iota_{{\cal X}_1}^*$ is onto, we have
$$\iota_{{\cal G}_1}^*(m)\cdot M_1=M_1\bullet \iota_{{\cal G}_1}^*(m),$$
for all $M_1\in LUC({\cal X}_1,{\cal G}_1)^*$. We now invoke Theorem \ref{5} to conclude that $\iota_{{\cal G}_1}^*(m)\in {\frak Z}({\cal X}_1,{\cal G}_1)$.}
\end{example}

\begin{example}\label{0925}
{\rm Let $Q_8$ be the quaternion group, the subgroup of the
general linear group $GL(2, \Bbb C )$ generated by the matrices
\begin{eqnarray*}
\bf{1}=\left(\begin{array}{clcr}
1 & 0\\
0 & 1
\end{array}
\right), \quad \bf{i}=\left(\begin{array}{clcr}
i & 0\\
0 & -i
\end{array}
\right), \quad \bf{j}=\left(\begin{array}{clcr}
0 & 1\\
-1 & 0
\end{array}
\right), \quad \bf{k}=\left(\begin{array}{clcr}
0 & i\\
i & 0
\end{array}
\right),
\end{eqnarray*}
which can be also presented by $$Q_8=\langle a,\,b:\,a^4=1,
\,a^2=b^2,\,b^{-1}ab=a^{-1}\rangle,$$
where $a=\bf{i}$ and $b=\bf{j}$ for
instance. Take ${\cal K}=\langle a\rangle$, which is a normal subgroup of $Q_8$.
If we set ${\cal G}=Q_8\times \Bbb F_2$, the direct product of $Q_8$ and
the free group on a two-element set, then ${\cal
X}={\cal K}\times \Bbb F_2$ is a closed normal subgroup of
discrete topological group $\cal G$ and so, $\cal G$ acts on ${\cal X}$ by $${\cal G}\times
{\cal X}\rightarrow {\cal X}:(s,{x})\mapsto
s{x}s^{-1}.$$ Since the algebraic centre ${\frak C}_{\cal G}({\cal G})$ contains
${\frak C}_{Q_8}(Q_8)=\{{\bf{i}}^2,\bf{1}\}$, by part (b) of Example~\ref{exmpl111}, we have $M({\cal
G})\varsubsetneq {\frak{Z}}({\cal X},{\cal G})$. On the other hand, if ${\frak{Z}}({\cal X},{\cal G})= LUC({\cal G})^*$, then by Example~\ref{rem1}, we would have
${\frak{Z}}({\Bbb F_2},{\Bbb F_2})=LUC(\Bbb F_2)^*$, which is impossible; This is because of, since $\Bbb F_2$ is not compact,
$${\frak{Z}}({\Bbb F_2},{\Bbb F_2})={\frak{Z}}(\Bbb F_2) \varsubsetneq LUC(\Bbb F_2)^*.$$
Hence, for $\cal G$-space $\cal X$, we have
$M({\cal G})\varsubsetneq {\frak{Z}}({\cal X},{\cal G})\varsubsetneq LUC({\cal G})^*$.

}
\end{example}

%%%%%%%%%%%%%%%%%%%%%%%%%%%%%%%%%%%%%%%%%%%%%%%%%%%%%%%%%%%%%%%%%%%%%%%%%%%%%%%%%%%%
%%%%%%%%%%%%%%%%%%%%%%%%%%%%%%%%%%%%%%%%%%%%%%%%%%%%%%%%%%%%%%%%%%%%%%%%%%%%%%%%%%%%

\vspace{3mm}

\noindent {\sc Hossein Javanshiri}\\
Department of Mathematics,
Yazd University,
P.O. Box: 89195-741, Yazd, Iran\\
E-mail: h.javanshiri@yazd.ac.ir\\

\noindent {\sc Narguess Tavallaei}\\
Department of Mathematics,
School of Mathematics and Computer Science,
Damghan University, Damghan, Iran.
E-mail: tavallaie@du.ac.ir

\end{document}